\documentclass{aupl}

\usepackage{
amsfonts,
amsmath,
latexsym,
amssymb,
enumerate,
verbatim,
mathrsfs,
centernot,
accents,
graphicx,
rotate,
color,
epigraph,
}

\newcommand{\labbel}[1]{\label{#1} [[{\bf #1}]]}  
\renewcommand{\labbel}{\label}

\newtheorem{theorem}{Theorem}

\newtheorem{proposition}[theorem]{Proposition}

\theoremstyle{definition}
\newtheorem{definition}[theorem]{Definition}

\theoremstyle{remark}
\newtheorem{remark}[theorem]{Remark}

\newtheorem{examples}[theorem]{Examples}

\renewcommand{\theequation}
{Eq. \arabic{equation}}

\begin{document}

\title{Exact-$m$-majority  terms}

\author{Paolo Lipparini} 
\address{Dipartimento 
 di Matematica\\Viale della  Ricerca  Esatta
\\Universit\`a di Roma ``Tor Vergata'' 
\\I-00133 ROME ITALY}
\urladdr{http://www.mat.uniroma2.it/\textasciitilde lipparin}


\thanks{Work performed under the auspices of G.N.S.A.G.A. Work 
partially supported by PRIN 2012 ``Logica, Modelli e Insiemi''.
The author acknowledges the MIUR Department Project awarded to the
Department of Mathematics, University of Rome Tor Vergata, CUP
E83C18000100006.
\\
\indent Date: \today}

\begin{abstract}
We say that an idempotent term $t$ is an 
\emph{exact-$m$-majority  term} 
if $t$ evaluates to $a$,
whenever the element $a$ occurs exactly 
$m$ times in the arguments of $t$, 
and all the other arguments are equal. 
     
If $m<n$ and some variety $\mathcal V$ has an
 $n$-ary  
exact-$m$-majority  term, then $\mathcal V$ is
congruence modular.
For certain values of $n$ and  $m$,
for example, $n=5$ and  $m=3$, the existence of 
an $n$-ary  
exact-$m$-majority  term 
neither implies congruence distributivity,
nor congruence permutability.
\end{abstract}

\maketitle

Near-unanimity terms
have been around in universal algebra
starting from  the 70's in the past century
\cite{BP,Mi} and recently have
played an important role
in tractability problems \cite{BIM,IMM}.
Recent results about near-unanimity terms
include  \cite{Bar,misharp,MZ}. 
Further references
 can be found in the quoted
works. 

A ternary near-unanimity term is a majority term.
Curiously enough, the opposite notion of a
minority term has proved quite interesting
\cite{KOVZ}. 
In the final section of \cite{misharp} 
we generalized 
the notion of a minority term
to  a ``lone-dissent'' term. 
A \emph{lone-dissent term} $u$ 
returns an element appearing just once
among its arguments, provided
all the other arguments are equal. 
In contrast, in the same situation, a near-unanimity term 
returns the element appearing all but one time.
In both cases, we are in the situation
in which the arguments of $u$
are always chosen from a pair of elements.
If $u$ returns the element appearing
$n-1$ times, where $u$ is $n$-ary,
then $u$ is near-unanimity term.
If $u$ returns the element appearing
once, then $u$ is a lone-dissent term. 
Thus the two notions have the following common generalization.

\begin{definition} \labbel{e}
Suppose that $0<m\leq n$. An $n$-ary
term $u$ is an \emph{exact-$m$-majority term}
(in some algebra or some variety) if the equations
\begin{equation}\labbel{ee}
u(x_1, x_2, \dots, x_n ) = x
  \end{equation} 
hold, whenever 
only the variables $x$
and  $y$ occur in  $\{ x_1, x_2, \dots, x_n  \} $ 
and the set $\{ \, i \leq n \mid x_i = x \,\} $
has exactly $  m$ elements.

We shall see that the two situations
$m < \frac{n}{2} $ and $m > \frac{n}{2} $ are rather 
different. In both cases congruence modularity
follows (if $m < n$). On the other hand, the existence of  an
exact-$m$-majority term implies 
congruence permutability for $m < \frac{n}{2} $,
while it
never implies 
congruence permutability, for
$m > \frac{n}{2} $.
\end{definition}   

\begin{remark} \labbel{rm}   
(a) Notice that if  \eqref{ee}
holds in some variety,
 we can take $x=y$
in \eqref{ee}, hence every 
exact-$m$-majority term is idempotent.
An $n$-ary  term $u$  
is idempotent if and only if 
it is an  exact-$n$-majority term, in the present terminology.

(b) It follows immediately from the above remark
that if $m < n$, then having an  exact-$m$-majority term is
 a nontrivial idempotent Maltsev condition
(compare \cite{T} and  \cite[Lemma 9.4(3)]{HMK}). 
As mentioned, we shall prove the stronger result that
the existence of 
an exact-$m$-majority term implies congruence modularity.

(c) If $u$ is  an $n$-ary term,
then  $u$ is a \emph{near-unanimity term}  
if and only if $u$ is an exact-$(n{-}1)$-majority term. 
A term  $u$ is a \emph{lone-dissent term}  \cite{misharp} 
if and only if $u$ is an exact-$1$-majority term.
In this terminology, a \emph{minority term} \cite{KOVZ}
is a $3$-ary   lone-dissent term.

(d) If $n$ is even, then
a variety $\mathcal V$ has an
 exact-$ \frac{n}{2} $-majority term if and only if $\mathcal V$ is trivial.

(e) Intuitively, the notion of a (non exact) $m$-majority
term could appear more natural; however we shall soon see
that this ``non exact'' notion provides nothing essentially  new,
since it turns out to be equivalent to the existence of a majority term---possibly,
of distinct arity.
If 
 $ \frac{n}{2}  <m\leq n$ we say that an $n$-ary
term $u$ is an \emph{$m$-majority term} if
$u(x_1, x_2, \dots, x_n ) = x$
holds, whenever 
$|\{ \, i \leq n \mid x_i = x \,\} | \geq m$.

If $m < n$ and $u$ is an $m$-majority term,
then $v(x_1, x_2, \dots, x_{m+1} )=
u(x_1, x_2, \allowbreak  \dots, x_m, x_{m+1}, x_{m+1}, \dots, x_{m+1})$ 
is an $m+1$-ary near-unanimity term, hence the notion of a
(non exact) $m$-majority
term seems to have little interests.
Let us point out, however, that the notion has been used, letting
$m$ vary, in order to construct the main counterexample in \cite{misharp}.
See  \cite[Section 2]{misharp}.
 \end{remark}

Let $(m,n)$ denote the greatest common divisor
of $m$ and $n$.  

\begin{proposition} \labbel{prop}
Suppose that $n \geq 3$ and
$0 < m < n$. 
  \begin{enumerate}   
 \item 
If $m < \frac{n}{2}$ and some variety $\mathcal V$ 
has an   $n$-ary  exact-$ m $-majority term, then $\mathcal V$
is congruence permutable. 
\item 
If $m > \frac{n}{2}$, then there is some variety
 $\mathcal V$ with an $n$-ary    exact-$ m $-majority term and
which
is not congruence permutable. 
 \item
If $k | (m,n)$ and some variety $\mathcal V$ 
has an  $n$-ary   exact-$ m $-majority term, then $\mathcal V$
has an  $ \frac{n}{k} $-ary   exact-$ \frac{m}{k}  $-majority term.
\item 
If $n- m $ divides $ n$ and some variety $\mathcal V$ 
has an  $n$-ary   exact-$ m $-majority term,
then $\mathcal V$ has an  $ \frac{n}{n-m} $-ary near-unanimity term,
in particular, $\mathcal V$ is congruence distributive.
\item
If $hm \equiv 1 \pmod q$ and  $n=m+ kq$,
for some $h,q,k \in \mathbb N$, 
then the term $u(x_1, x_2, \dots, x_n ) = hx_1+ hx_2 + \dots + hx_n$ is
an  $n$-ary   exact-$ m $-majority term in any 
 abelian group of exponent dividing $q$ (in additive notation).
In particular, for such values of $m$  and $n$, if $q>1$ the existence of
an  $n$-ary   exact-$ m $-majority term  
does not imply congruence distributivity.
  \end{enumerate} 
 \end{proposition}

  \begin{proof}
(1) If $m < \frac{n}{2}$ and $u$ is  
an   $n$-ary  exact-$ m $-majority term, consider the term
 $t(x,y,z)= u(x, \dots , x, y, \dots, y, z, \dots, z)$,
where  both $x$ and  $z$ appear $m$ times and  
$y$ appears  $n-2m$ times. 
Notice that $n-2m>0$, since $m < \frac{n}{2}$.
Then $t$  is a Maltsev term witnessing
congruence permutability.

(2)  If $m > \frac{n}{2}$, then in lattices the term 
\begin{equation*}
u _{m,n}(x_1, \dots, x_n)= \prod _{|J|=m} \sum _{i \in J} x_i 
   \end{equation*}    
($J$ varying on  subsets of 
$\{ 1, \dots, n \})$
is an 
exact-$ m $-majority term, actually,
an $ m $-majority term.
However, lattices are not congruence permutable.

(3)  If $u$ witnesses the assumptions and
$\ell= \frac{n}{k}$, let
 $t(x_1, x_2, \dots, x_\ell  )
= u(x_1,x_1, \allowbreak \dots, x_2,x_2,\dots, x_\ell,x_\ell, \dots)$,
where each variable appears $k$ times on the
right-hand side.  

(4) Take $k= n-m $ in (3), then notice that 
$\frac{m}{k} + 1 = \frac{n}{k}$. 

(5) If $x$ appears $m$ times in the arguments of the term,
the sum of the corresponding summands gives
  $hmx$, that is, $x$, since 
 $hm \equiv 1 \pmod q$ and the group
is abelian  of exponent dividing $q$.
If the only other variable is $y$, it occurs $kq$ times,
hence the outcome is $0$, again  
  since the group
is abelian  and its  exponent divides $q$.
 \end{proof}

\begin{examples} \labbel{ex}
(a) The existence of a $5$-ary exact-$3$-majority term
does not imply congruence permutability, by Proposition \ref{prop}(2). 
It does not imply congruence distributivity, either, by Proposition \ref{prop}(5),
taking $h=k=1$ and $q=2$.  

Henceforth, if $\mathcal V_{3,5}$ 
is the variety with a $5$-ary operation satisfying 
the equations for a  $5$-ary exact-$3$-majority term,
then $\mathcal V_{3,5}$ is neither congruence permutable nor
congruence distributive. On the other hand, we shall show that every variety 
with an  exact-$m$-majority term is congruence modular. Hence
$\mathcal V_{3,5}$ seems to be an interesting example of 
a congruence modular variety which is neither  congruence permutable nor
congruence distributive.

(b)
If $m \neq \frac{n}{2}$, then there is a nontrivial algebra
(hence a nontrivial variety) with an 
$n$-ary exact-$m$-majority term.
Actually, every set admits an $n$-ary exact-$m$-majority operation.

Let $A$ be any nonempty set and fix some $a \in A$. Define
\begin{equation*}\labbel{ca}
u(a_1, \dots, a_n)     =\begin{cases} b &    \text{if  
$|  \{ \, i \leq n \mid a_i=b\, \}| = m  $
and $|\{ \, a_i \mid i \leq n  \,\}| =2$},
\\ b &    \text{if  
$a_1=a_2=\dots= a_n=b$},
\\  a & \text{otherwise.}
 \end{cases}
 \end{equation*} 

Since  $m \neq \frac{n}{2}$, the first clause provides a good definition.

(c) There is no nontrivial group $\mathbf G$
 with a $6$-ary exact-$2$-majority term.
To prove this, let us use multiplicative notation. 
A group term is a product of variables raised to some power;
if some term is evaluated for just one element $g$ and for the identity,
the outcome of the term is $g$ raised to the sum of the powers of
the occurrences of the corresponding variables.
If $t$ is a $6$-ary exact-$2$-majority term, then $t(g,g,e,e,e,e)=g$
and   $t(e,e,g,g,e,e)=g$, for every element $g \in G$.
By the preceding comment, if $h$ is the sum of the exponents
of all the occurrences of the first two variables, then $g^h=g$
in $\mathbf G$. Similarly $g^k=g$, where $k$ is the sum of the exponents
of all the occurrences of the third and fourth variables.     
 Hence $t(g,g,g,g,e,e) $ evaluates to $g ^{h+k}= g^hg^k=  g^2$ in $\mathbf G$.
 But also 
$t(g,g,g,g,e,e)=e$, by the majority assumption, thus 
$\mathbf G$ has exponent $2$, since the above
argument applies to every $g \in G$. Since every group of
exponent $2$ is abelian, every term of $\mathbf G$
can be represented as a product of variables.
 Since $\mathbf G$ is nontrivial 
of exponent $2$ and $t(g,g,e,e,e,e)=g$, for every $g \in G$,
then either $t$ does not depend on the first variable, or
 $t$ does not depend on the second variable. 
 Similarly, since   $t(e,e,g,g,e,e)=g$,  
then either $t$ does not depend on the third variable, or
 $t$ does not depend on the fourth variable.
Say, $t$ does not depend on the first and on the third variables.
Since $\mathbf G$ is non trivial, then there is  $g \in G$ with
$g \neq e$, but then we  have 
$t(g,e,g,e,e,e)=e \neq g$, a contradiction.

Notice that we have used only $4$ instances of the 
 majority rule (among 15 total instances). 

(d) By Proposition \ref{prop} (3), every variety 
 with a $6$-ary exact-$2$-majority term has 
 a $3$-ary exact-$1$-majority term, i.e.,
a minority term. In a group of exponent $2$
the term $xyz$ is a minority term.
On the other hand, by the previous item,
no nontrivial group has a  $6$-ary exact-$2$-majority term.
Thus there is a variety 
 with a 
 $3$-ary exact-$1$-majority term but without a
$6$-ary exact-$2$-majority term.
\end{examples}

\begin{theorem} \labbel{thm}
Suppose that $n \geq 3$ and
$0 < m < n$. Then every variety with  an
$n$-ary exact-$m$-majority term is congruence modular.
 \end{theorem} 

\begin{proof} 
If  $m < \frac{n}{2}$, then $\mathcal V$ is congruence permutable,
by Proposition \ref{prop}(1), hence $\mathcal V$ is congruence modular. 
If $m = \frac{n}{2}$, then $\mathcal V$ is a trivial variety, by Remark \ref{rm}(d).

It remains to deal with the case  $m > \frac{n}{2}$, so let
$u$ be an  $n$-ary exact-$m$-majority term in this case.
Let $k=n-m$ and $h$ be the remainder
of the division of $n$ by $k$
(the proof below works also in case $h=0$; anyway,
the case $h=0$ is already covered by Proposition \ref{prop}(4)).
  Let $x^k$ be an abbreviation for the expression
``$x,x, \dots, x$'', with $k$ occurrences of $x$. Consider the terms
\begin{align*} 
d_1(x,y,z)&= u( x^h,x^k,x^k,x^k,\dots, x^k,x^k,x^k,y^k,z^k)
\\
d_2(x,y,z)&= u( x^h,x^k,x^k,x^k,\dots, x^k,x^k,y^k,z^k,z^k)
\\
d_3(x,y,z)&= u( x^h,x^k,x^k,x^k,\dots, x^k,y^k,z^k,z^k,z^k)
\\
&\dots
\\
d_{\ell-2}(x,y,z)&= u( x^h,x^k,y^k,z^k,\dots, z^k,z^k,z^k,z^k,z^k)
\\
d_{\ell-1}(x,y,z)&= u( x^h,y^k,z^k,z^k,\dots, z^k,z^k,z^k,z^k,z^k)
\\
q(x,y,z)&= u( x^h,y^{k-h},z^h,z^k,z^k,\dots, z^k,z^k,z^k,z^k,z^k)
 \end{align*}
were $\ell$ is the integer quotient of the division of $n$ by $k$.
Notice that $\ell >1$, since  $m > \frac{n}{2}$.
The above terms satisfy   $d_i(x,z,z)=d_{i+1}(x,x,z)$, for $1 \leq i < \ell -1$.
If $\ell=2$, the terms to be considered are
\begin{align*} 
d_{1}(x,y,z)&= u( x^h,y^k,z^k)
\\
q(x,y,z)&= u( x^h,y^{k-h},z^h,z^k)
 \end{align*}
All the above terms satisfy    $d_{\ell-1}(x,z,z)=q(x,z,z)$. 
Due to the exact $m $ ($ =n-k$) majority rule, the above terms also satisfy
$d_i(x,y,x)=x$, for $1 \leq i \leq \ell -1$,  $x=d_1(x,x,z)$  and
$q(x,x,z)=z$. In the terminology from \cite[p. 205]{KKMM}, the terms 
  $d_1$, \dots $d_{\ell-1}$, $q$ are \emph{directed Gumm terms},
and the existence of such a sequence of terms implies congruence modularity,
by the easy part of \cite[Theorem 1.1, Clause 3]{KKMM}. 
\end{proof}  

\begin{remark} \labbel{edge}
If $m < n$, an $n$-ary exact-$m$-majority term
is a $\Delta$-special cube term, in the terminology of \cite[Definition 2.7]{BIM},
for some appropriate $\Delta$.
Indeed, if we write down all the equations defining  
an $n$-ary exact-$m$-majority term, we get a matrix 
with $k={n \choose m}$ rows, and only the last column is
constantly $x$.
Hence a variety with an $n$-ary exact-$m$-majority term
has a $k$-edge term, by \cite[Theorem 2.12]{BIM}.
 In particular, \cite[Theorem 4.2]{BIM} furnishes another proof of congruence modularity.
 \end{remark}

\newcommand{\eeeemph}{\emph}


\begin{thebibliography}{10}


\bibitem{BP}
 K. A. Baker, A. F. Pixley,  \emph{Polynomial interpolation and the Chinese remainder
theorem for algebraic systems}, Math. Z. \textbf{143}, 165--174 (1975) 

\bibitem{Bar} L. Barto, \emph{Finitely related algebras in congruence distributive varieties have near unanimity terms}, Canad. J. Math. \textbf{65}, 3--21 (2013) 


\bibitem{BIM} J. Berman, P.  Idziak, P. Markovi\'c, R. McKenzie, M. Valeriote, 
 R. Willard, 
\emph{Varieties with few subalgebras of powers}, Trans. Amer. Math. Soc. \textbf{362},
1445--1473 (2010) 


\bibitem
{HMK}
D. Hobby, R. McKenzie, \emph{The structure of finite algebras},
Contemp. Math.  \textbf{76} (1988)

\bibitem{IMM}
P.  Idziak, P. Markovi\'c, R. McKenzie, M. Valeriote, 
 R. Willard, 
\eeeemph{Tractability and learnability arising from algebras with few
 subpowers},
{SIAM J. Comput.}
\textbf{39},
3023--3037 (2010)



\bibitem{KKMM}
A. Kazda, M.  Kozik, R.  McKenzie, M. Moore,
 \eeeemph{Absorption and directed J\'onsson terms},
in: J.\ Czelakowski  (ed.), \eeeemph{Don Pigozzi on Abstract Algebraic Logic, Universal Algebra, and Computer Science}, Outstanding Contributions to Logic \textbf{16}, Springer, Cham, 203--220  (2018)

\bibitem{KOVZ} 
 A. Kazda,  J. Opr\v sal,  M. Valeriote, D. Zhuk, 
\emph{Deciding the existence of minority terms},
Canad. Math.
Bull. \textbf{63}, 577--591 (2020)



\bibitem{misharp}
P. Lipparini,  \emph{Mitschke's Theorem is sharp},
 {Algebra Universalis}  \textbf{83}, 7 (2022) 1--20


 
\bibitem{MZ}
M. Mar\'{o}ti, L. Z\'{a}dori, 
\eeeemph{Reflexive digraphs with near unanimity polymorphisms},
 {Discrete Math.}
\textbf{312},
2316--2328 (2012)
 

\bibitem{Mi} 
A. Mitschke, 
\emph{Near unanimity identities and congruence distributivity in equational classes},
 Algebra Universalis \textbf{8},  29--32  (1978)



\bibitem{T} W. Taylor, 
\emph{Varieties obeying homotopy laws}, Canadian J. Math. 29 (1977),
498-527




\end{thebibliography}
\end{document}